%Paper: alg-geom/9301005
%From: <GUES%UORDB5.bitnet@CUNYVM.CUNY.EDU>
%Date: Sun, 24 Jan 93 14:27 EDT
%Date (revised): Wed, 24 Aug 94 08:06 EDT

%%%%%%%%%%%%%%%%%%%%%%%%%
%%%%%% AMS-TeX 2.1 %%%%%%
%%%%%%%%%%%%%%%%%%%%%%%%%

\input amstex
\input amsppt.sty

\refstyle{A}

\magnification=1200
\parskip 10pt
\pagewidth{5.4in}
\NoBlackBoxes
\leftheadtext{Rational Curves on a Toric Variety}
\rightheadtext{Rational Curves on a Toric Variety}

\redefine\Bbb{\bold}
\define\C{\Bbb C}
\define\R{\Bbb R}
\define\Z{\Bbb Z}

\define\CP{\C P}
\define\TC{{T^{\C}}}
\define\CS{{\C^\ast}}
\define\QXD{Q^X_D}

\define\hX{\hat X}

\define\hD{\hat D}
\define\hDe{\hat\De}

\define\al{\alpha}
\define\be{\beta}
\define\de{\delta}

\define\De{\Delta}

\define\si{\sigma}
\define\ta{\tau}
\define\th{\theta}

\define\sub{\subseteq}

\define\es{\emptyset}
\redefine\ll{\lq\lq}
\redefine\rr{\rq\rq\ }
\define\rrr{\rq\rq}

\define\lan{\langle}
\define\ran{\rangle}
\define\pr{\prime}

\redefine\deg{\operatorname {deg}}
\define\Hol{\operatorname {Hol}}
\define\Map{\operatorname {Map}}

\define\Hom{\operatorname {Hom}}

\define\Ker{\operatorname{Ker}}
\redefine\min{\operatorname{min}}
\redefine\max{\operatorname{max}}
\define\SF{\operatorname {SF}}
\define\Pic{\operatorname {Pic}}
\define\Gr{\operatorname {Gr}}

\topmatter
\title The topology of the space of\\
rational curves on a toric variety
\endtitle
\author Martin A. Guest
\endauthor
\address Department of Mathematics, University of
Rochester, Rochester, New York 14627, USA
\endaddress
\email gues\@db1.cc.rochester.edu
\endemail
\endtopmatter

\document

\baselineskip=13pt

\subheading{\S 1. Introduction}

A toric variety is, roughly speaking, a complex algebraic variety
which is
the (partial) compactification of an algebraic torus $\TC=(\CS)^r$. It
admits (by definition)
an action of $\TC$ such that, for some point $\ast\in X$, the orbit of
$\ast$
is an embedded copy of $\TC$. The most significant property of a
toric
variety is the fact that it is characterized entirely by a combinatorial
object, namely its {\it fan}, which is a collection of convex cones in
$\R^r$.   As a general
reference for the theory of toric varieties we use
\cite{Od1}, together with the recent lecture notes \cite{Fu}.

In this paper we shall study the space of rational curves on a compact
 toric variety $X$. We shall obtain a configuration space description
of the space $\Hol(S^2,X)$ of all holomorphic
(equivalently, algebraic) maps from the Riemann sphere
$S^2=\C\cup\infty$ to $X$.  Our main application of this concerns
fixed components $\Hol^\ast_D(S^2,X)$  of
$\Hol^\ast(S^2,X)$, where the symbol $D$ will be
explained later, and where the asterisk indicates that the maps are
required to satisfy the condition $f(\infty)=\ast$.
If $\Map^\ast_D(S^2,X)$ denotes the corresponding space of
continuous maps, we shall show that
{\it the inclusion
$$
\Hol^\ast_D(S^2,X)\to\Map^\ast_D(S^2,X)
$$
induces isomorphisms of homotopy groups up to some dimension
$n(D)$}, and we shall give a procedure for computing $n(D)$.

A theorem of
this type was proved in the case $X=\CP^n$ by Segal (\cite{Se}), and
indeed that theorem provided the motivation for the present work. Our
main idea is that the result of Segal may be interpreted as a result
about configurations of distinct points in $\C$ which have labels in a
certain {\it partial monoid}. We shall show that $\Hol_D^\ast(S^2,X)$
may be identified with a space $Q_D^X(\C)$ of configurations of
distinct points in $\C$ which have labels in a partial monoid $M_X$,
where $M_X$ is derived from the fan of $X$. Then we shall extend
Segal's method so that it applies to this situation.

A feature of the method  is the idea that the functor $U\mapsto\pi_i
Q_D^X(U)$ resembles a homology theory.  This
functor has some similarities with the Lawson homology functor
introduced in \cite{La}, in the sense that both are generalizations of
the Dold-Thom functor $U\mapsto\pi_i Sp^d(U)$, where $Sp^d(U)$ is
the $d$-th symmetric product of $U$. This latter space can be
considered as (a subspace of) the space of
configurations of distinct points in $U$ which have labels in the
partial  monoid $\{1,2,\dots,d\}$.
(The label of a point, here, is simply its
multiplicity.)  The Dold-Thom functor resembles the ordinary
homology functor, in the sense that
$\pi_iSp^d(U)\cong \tilde H_iU$, for $i$
less than some dimension $n(d,U)$  (\cite{DT}).
It is important to note that the labelled configurations here are
topologized so that when two distinct points \ll collide\rrr, their
labels are added; if the addition is not defined, then the collision is
prohibited.

A second aspect concerns a well known problem inspired by Morse
theory. To explain this, we note that the above theorem
(with $n(D)$ non-trivial!) is definitely not valid for arbitrary compact
complex analytic spaces, or even complex manifolds. For example,
there are no non-constant holomorphic maps from $S^2$ to the Hopf
surface $S^1\times S^3$, or indeed to any Abelian variety.
Nevertheless, there is reason to believe that a theorem of the above
type might hold for compact K\"ahler manifolds, because in this case
the holomorphic maps (in a fixed component of smooth maps) are
precisely the absolute minima of the energy functional.  A suitable
extension of Morse theory would then explain such a result (although
there would be no guarantee that $n(D)$ would be non-zero; for
example it is known that there exist K\"ahler manifolds with very
few rational curves).  Our results confirm this Morse theory principle
for smooth toric varieties, and they provide some evidence that it
extends even to certain {\it singular} varieties.

This paper is arranged as follows.
After a brief review of toric varieties in \S 2, we proceed to describe
the correspondence between holomorphic maps and labelled
configurations in \S 3, in the case of a projective toric variety
(singular or not). The proof of the main theorem
in the case of a
non-singular projective toric variety (Theorem 4.1) is given in \S 4.
It falls into three parts. First, we show that the
homotopy groups of $Q^X_D(\C)$ \ll stabilize\rr as $D\to\infty$. This
can be
reduced to the corresponding fact for the symmetric product. The
method we
use here is based on \cite{GKY},\cite{Gu1} as the method used in
\cite{Se} for the case $X=\CP^n$ does not seem to extend to the case
of general $X$.  Second, we show
(using the homology-like properties of $\pi_i Q^X_D(\C)$)
that $Q^X_\infty(\C)$ is actually
homotopy equivalent to a component of $\Map(S^2,X)$. This idea, due
originally to Gromov and Segal, has been used several times in the
past, e.g.  in \cite{Mc},\cite{Se},\cite{Gu1},\cite{Gu2}. Third, we show
that this homotopy  equivalence actually arises from the inclusion
map of the theorem.
In \S 5 we discuss the considerably more
difficult case of singular projective toric varieties.  We are not able
to give a general result here which covers all cases, but we shall
establish the result in a basic situation (Theorems 5.1 and 5.2) and
then illustrate by examples how our method can in principle be
applied to the remaining cases.  We also sketch how the
results may be generalized
to non-projective toric varieties. For technical reasons we study only
compact toric varieties in this paper, although it seems likely that a
similar method also works in the non-compact case (see \cite{GKY}
for an example).

The author is indebted to A. Kozlowski and K. Yamaguchi for numerous
helpful conversations concerning configuration spaces. He is also very
grateful to D. Cox, M. McConnell, T. Oda, M. Oka, B. Sturmfels and A.
Vitter for their comments.  In addition, Professor Oda suggested
a number of improvements to earlier versions of this paper, and
these were greatly appreciated.

\subheading{\S 2. Toric varieties}

We shall summarize some of the basic properties of toric varieties,
from
\cite{Od1}. Let $X$ be an irreducible normal algebraic variety. One
says that
$X$ is a {\it toric variety} if it has an algebraic action of an algebraic
torus
$\TC=(\CS)^r$, such that the orbit $\TC\cdot\ast$ of some point
$\ast\in X$
is dense in $X$ and isomorphic to $\TC$.

A toric variety is characterized up to isomorphism by its {\it fan},
which is
a finite collection $\Delta$ of strongly convex rational polyhedral
cones in $\R^r$ such that every face of an element of $\Delta$ belongs
to $\De$, and the intersection of any two elements of $\Delta$ is a
face of each. (A strongly convex rational polyhedral cone in $\R^r$ is a
subset of $\R^r$ of the form $\{\sum_{i=1}^sa_in_i\ \vert\ a_i\ge
0\}$, where $\{n_1,\dots,n_s\}\sub\Z^r$, which does not contain any
line.)
Given a fan $\Delta$, an associated toric
variety may be constructed abstractly as a union of
affine varieties $U_\si$,  $\si\in\De$.

It is possible to give a concrete description of {\it projective} toric
varieties, as follows. Let  $m_1,\dots,m_N\in\Z^r$, such that the
elements $m_i-m_j$ generate $\Z^r$.  Consider the action of
$(\C^\ast)^r$ on $\CP^N$ given by the formula
$u\cdot[z_0;\dots;z_N]=[z_0;u^{m_1}z_1;\dots;u^{m_N}z_N]$, where
$u=(u_1,\dots,u_r)$, $m_i=((m_i)_1,\dots,(m_i)_r)$,  and
$u^{m_i}=u_1^{(m_i)_1}\dots u_r^{(m_i)_r}$. Then the closure of the
orbit of $[1;\dots;1]$ is a toric variety.  This gives rise to a
second
characterization of a toric variety embedded in projective space,
namely that it is defined by equations of the form \ll monomial in
$z_0,\dots,z_N$ = monomial in $z_0,\dots,z_N$\rrr. A third explicit
description will be mentioned at the end of \S 5.

{}From the general construction, it follows that there is a one to
one correspondence between $\TC$-orbits of codimension $i$ in $X$
and cones of dimension $i$ in $\De$. The closure of the $\TC$-orbit
corresponding to a cone $\si$ is a toric subvariety; it is the union of
the  orbits corresponding to those $\ta\in\De$ such that $\si$ is a
face of $\ta$.
Not surprisingly, geometrical and topological properties of $X$ are
reflected
in the fan $\De$. For example, $X$ is non-singular if and only if, for
each $\si\in\De$, the generators $n_1,\dots,n_s$ can be extended to a
generating set for $\Z^r$ (\cite{Od1}, Theorem 1.10).  The variety $X$
is compact if and
only if $\bigcup_{\si\in\De}\si=\R^r$ (\cite {Od1}, Theorem 1.11).
{\it From now on we shall assume that $X$ is a compact toric
variety.}

We shall be particularly concerned with the topology of $X$, and
with line bundles over $X$.
It is known that the fundamental group of any toric
variety $X$ is isomorphic to the quotient of $\Z^r$ by the subgroup
generated by $\bigcup_{\si\in\De}\si\cap\Z^r$ (\cite{Od1},
Proposition 1.9).
{}From this and the compactness criterion, $\pi_1X=0$.
To get further information, we need to introduce some more
notation.  Let $\si_1,\dots,\si_u$ be the
one-dimensional cones in $\De$.  We have $\si_i\cap\Z^r=\Z v_i$ for
some $v_i\in\Z^r$.  Let $X_1,\dots,X_u$ be the closures of the
corresponding (codimension one) $\TC$-orbits in $X$.
Equivariant line bundles on $X$ correspond to invariant  Cartier
divisors on $X$.  If $X$ is compact, these correspond to
\ll $\De$-linear support functions\rrr, i.e. functions $h:\De\to\R$
which are  linear on each cone $\si$ and $\Z$-valued on $\De\cap\Z^r$
(\cite{Od1},  Proposition 2.4).  Let $\SF(\De)$ denote the group of
$\De$-linear support functions. For $h\in \SF(\De)$, a divisor of the
corresponding line bundle is given by  $\sum_{i=1}^u h(v_i)X_i$.
Hence we obtain an inclusion $\SF(\De)\to\bigoplus_{i=1}^u\Z\si_i$,
$h\mapsto\sum_{i=1}^u h(v_i)\si_i$.  From now on, we shall {\it
identify} $\SF(\De)$ with a subgroup of $\bigoplus_{i=1}^u\Z\si_i$. It
represents the subgroup consisting of invariant
Cartier divisors of the group
$\bigoplus_{i=1}^u\Z\si_i$ of invariant Weyl divisors.  The
inclusion is an isomorphism if $X$ is non-singular (\cite{Od1},
Proposition 2.1).

We have another natural inclusion $\iota:\Z^r\to \SF(\De)$,
$m\mapsto\langle m,\ \ \rangle$,  and the quotient group  is
isomorphic to the Picard group $\Pic(X)\cong H^1(X,\Cal O_X^\ast)$, if
$X$  is compact (\cite{Od1}, Corollary 2.5).  Moreover, in this case, it
is known that $H^i(X,\Cal O_X)=0$ for $i\ge 1$ (\cite {Od1}, Corollary
2.8), so we have $H^2X\cong H^1(X,\Cal O_X^\ast)\cong \Pic(X)\cong
\SF(\De)/\Z^r$.
Since we consider maps $S^2\to X$ in this paper, we shall need a
description of the group $\pi_2X$. If
$H_2X$ is torsion free, then $H_2X\cong(H^2X)^\ast$, and so we have
$$
\pi_2X\cong H_2X\cong
(H^2X)^\ast\cong(\SF(\De)/\Z^r)^\ast\cong\Ker\iota^\ast
$$
where $\iota^\ast:\SF(\De)^\ast\to(\Z^r)^\ast$ is the dual of $\iota$.
If $X$ is non-singular then $H_2X$ is torsion free, because
the integral cohomology of $X$ is torsion free, by a theorem of
Jurkiewicz-Danilov (\cite{Od1}, page 134).

\noindent {\it Example 2.1: Complex projective space, $X=\CP^n$.}
Let $E_1,\dots,E_n$ be the standard orthonormal basis of
$\R^n$.
Let $\rho_0,\dots,\rho_n$ be the one-dimensional cones spanned by
$E_0=-\sum^n_{i=1}E_i,E_1,\dots,E_n$ (respectively). We obtain a fan
in
$\R^n$ by taking the cones spanned by all proper subsets of
$\{\rho_0,\dots,\rho_n\}$ (together with the zero-dimensional cone
given
by the origin). The associated toric variety is isomorphic to $\CP^n$.
In terms of the usual homogeneous coordinates for $\C P^n$,
the algebraic torus $\TC=(\CS)^n$ acts by the formula
$(v_1,\dots,v_n)\cdot [z_0;\dots;z_n]=[z_0;v_1z_1;\dots;v_nz_n]$.
Let $\ast=[1;\dots;1]$. The closures of the codimension
one $\TC$-orbits are the hyperplanes $P_0,\dots,P_n$, where $P_i$ is
defined by the condition $z_i=0$.

\noindent {\it Example 2.2: The Hirzebruch surface $X=\Sigma_k$
(\cite{Od1}, page 9, Example (iii)).}
Consider the fan in $\R^2$ given by the four
two-dimensional cones (and all their faces) spanned by the four
vectors
$v_1=(1,0)$,
$v_2=(0,1)$,
$v_3=(-1,k)$,
$v_4=(0,-1)$.
Thus, the one-dimensional cones in this fan are
$\si_i=\R_{\ge0}v_i$  for $i=1,2,3,4$,
 and we have $\si_i\cap\Z^2=\Z v_i$.
The construction produces a variety isomorphic to
$\Sigma_k$.
The classes $[\si_1],[\si_2],[\si_3],[\si_4]$ in
$\bigoplus_{i=1}^4\Z\si_i/\Z^2$ satisfy the
relations $[\si_1]=[\si_3]$, $[\si_4]=[\si_2]+k[\si_3]$ (corresponding
to the generators $(1,0),(0,1)$ of $\Z^2$). Hence
$H^2\Sigma_k\cong\Z\oplus\Z$.
To describe this variety directly, let us consider the embedding
$$
\Sigma_k=\{([x_0;x_1;x_2],[y_1;y_2])\ \vert\ x_1y_1^k=x_2y_2^k\}
\sub \C P^1\times\C P^2.
$$
The algebraic torus $\TC=(\CS)^2$ acts on $\Sigma_k$ as follows:
$$
(v_1,v_2)\cdot([x_0;x_1;x_2],[y_1;y_2])=
([v_1x_0;v_2^kx_1;x_2],[y_1;v_2y_2]).
$$
Let $\ast=([1;1;1],[1;1])$. The closures of the
codimension one $\TC$-orbits are the
four embedded copies of $\CP^1$ defined by $X_1=\{x_2=0,y_1=0\}$,
$X_2=\{x_1=0,x_2=0\}$, $X_3=\{x_1=0, y_2=0\}$, $X_4=\{x_0=0\}$.
The natural projection $\Sigma_k\to\CP^1$ exhibits $\Sigma_k$ as
$\Bbb P(\Cal O(0)\oplus\Cal O(-k))$, which is the bundle obtained from
$\Cal O(-k)$ by fibre-wise one point compactification. The
$0$-section and $\infty$-section are given by $X_2$ and $X_4$, and
the fibres over $[0;1],[1;0]$ are given by $X_1,X_3$.

\noindent{\it Example 2.3: The weighted projective spaces
$X=P(a_0,\dots,a_n)$ (\cite{Fu}, Section 2.2).}
The weighted projective space $P(a_0,\dots,a_n)$ is
defined to be the quotient of $\CP^n$ by the action of the finite group
$(\Z/a_0\Z)\times\dots\times(\Z/a_n\Z)$ given by
$$
(\omega_0,\dots,\omega_n).[z_0;\dots;z_n]=[\omega_0z_0;\dots;
\omega_nz_n],
$$
where $\omega_i$ is a primitive $a_i$-th root of unity. Without loss
of generality we may assume $a_0=1$. In this case, a suitable fan is
generated by the vectors
$\ -\sum_{i=1}^n a_i E_i,E_1,\dots,E_n$, and we have
$H^2P(a_0,\dots,a_n)\cong\Z$. These varieties may have
singularities.

\noindent{\it Example 2.4: Compact non-singular toric surfaces.}
These are classified in \cite{Fu}, Section 2.6, and \cite{Od1},
Sections 1.6,1.7. They are obtained from $\C P^2$ or $\Sigma_k$ by
blowing up a finite number of fixed points of the torus
action.

\noindent{\it Example 2.5: The closure of an algebraic torus orbit in
a (generalized) flag manifold.}
A Lie-theoretic description of the fan is given in \cite{FH};
see also \cite{Da} and \cite{Od2}, and the references given in
\cite{Od2} to closely related work of the Gelfand school.  It was
pointed out in \cite{Od2} that the normality of these varieties
remains to be verified. This omission has recently been rectified in
\cite{Da}, in the case of generic orbits.
These varieties may have singularities.
As a concrete example, we mention the famous \ll tetrahedral
complex\rrr, which is a singular three-dimensional subvariety $X$ of
the Grassmannian $\Gr_2(\C^4)$.  For historical remarks, including a
description of the role played by this variety in the early development
of Lie theory, we refer to Section 4.2 of \cite{GM}. In mundane terms,
if $\Gr_2(\C^4)$ is realized as the subvariety of $\C P^5$ given by the
usual Pl\"ucker equation $z_0z_1-z_2z_3+z_4z_5=0$, then $X$ is
given by the equations $z_0z_1=\al z_2z_3=\be z_4z_5$, where
$\al,\be$ are fixed complex numbers such that
$1-\al^{-1}+\be^{-1}=0$ and $\al,\be\ne0,1,\infty$. If, on the other
hand, $\Gr_2(\C^4)$ is considered as a generalized flag manifold of
the group $SL_4(\C)$, and if $(\C^\ast)^3$ is considered in the usual
way to be a maximal torus of $SL_4(\C)$, which therefore acts
naturally on $\Gr_2(\C^4)$, then $X$ occurs as the closure of a
generic orbit. It follows that the fan of $X$ can be obtained from the
results of \S 4 of \cite{FH}. After some re-normalization, it is the
fan in $\R^3$ with six
three-dimensional cones (and all their faces) spanned by the vectors
$(\pm1,\pm1,\pm1)$. The lattice is taken to be that which is
generated by $(\pm1,\pm1,\pm1)$, however, rather than $\Z^3$. One
has $H^2X\cong\Z$.

\subheading{\S 3. The configuration space for projective toric
varieties}

As in \S 2, let $\si_1,\dots,\si_u$ be the one-dimensional cones in
the fan  $\De$, and let $X_1,\dots,X_u$ be the closures of the
codimension one $\TC$-orbits in $X$.  Thus, $X_1\cup\dots\cup X_u$
is the complement of the \ll big orbit\rr $\TC\cdot\ast$.  We shall
assume as usual that $H_2X$ is torsion free, in order to make use of
the description of $\pi_2X$ which was given in $\S 2$.

If $f$ is a holomorphic map such that $f(\infty)=\ast$, then $f(S^2)$
is not contained in any of the subvarieties $X_i$, and so $f(S^2)\cap
X_i$ must be a finite (possibly empty) set of
points.  We associate to $f$ the finite set of distinct points
$\{z\in\C\ \vert\ f(z)\notin\TC\cdot\ast\}$, and to each such point
$z$ we associate --- provisionally --- a \ll label\rr
$l_z=((l_z)_1,\dots,(l_z)_u)$, where $(l_z)_i$ is the non-negative
integer given by the (suitably defined) intersection number of $f$ and
$X_i$ at $z$.

It turns out to be more natural to regard the label $l_z$ as an element
of $\SF(\De)^\ast$, i.e. $\Hom(\SF(\De),\Z)$. Therefore, the
provisional definition of the labelled configuration associated to $f$
will be replaced by the following construction. Let
$Q(\C;\SF(\De)^\ast)$ be the space of configurations of distinct
points in $\C$, where the points have labels in the group
$\SF(\De)^\ast$. An element of $Q(\C;SF(\De)^\ast)$ may be written
in the form $\{(z,l_z)\}_{z\in I}$, where $I$ is a
finite subset of $\C$ and
$\{l_z\}_{z\in I}\sub SF(\De)^\ast$.
There is a natural topology on $Q(\C;\SF(\De)^\ast)$ which
permits two distinct points in a configuration to \ll coalesce\rrr,
whereupon their labels are added.   Thus, $Q(\C;\SF(\De)^\ast)$
consists of a collection of disconnected components
$Q_D(\C;\SF(\De)^\ast)$, indexed by elements $D=\sum_{z\in I}l_z$ of
$SF(\Delta)^\ast$. Each component is contractible as all particles
may be moved to the origin.

\noindent{\bf Definition:} To a holomorphic map
$f\in\Hol^\ast(S^2,X)$ we associate a configuration in
$Q(\C;\SF(\De)^\ast)$ by means of the map
$$
\al^X:\Hol^\ast(S^2,X) \longrightarrow Q(\C;\SF(\De)^\ast),\quad
\al^X(f)=\{(z,l_z)\ \vert\ f(z)\notin T^{\C}\cdot\ast\}
$$
where, for any (Cartier) divisor $\tau\in \SF(\De)$,
$l_z(\tau)$ is the multiplicity of $z$ in the divisor
$f^{-1}(\tau)$.

\noindent Our main observation will be that the map $\al^X$ is a
homeomorphism to its image, and that the image has a simple
characterization. To obtain this characterization, we observe that
the configuration obtained from a map $f$ must satisfy two kinds of
properties.

First, geometry forces the following conditions on the label $l_z$ of
a point $z$:

\noindent (X)\qquad If $\tau\ge0$, then $l_z(\tau)\ge0$.
If $\tau_{i_1}\cap\dots\cap\tau_{i_j}=\es$, then
$l_z(\tau_{i_1})\dots l_z(\tau_{i_j})=0$
(i.e. $l_z(\tau_{i_1}),\dots, l_z(\tau_{i_j})$ cannot all be non-zero).

\noindent If $X$ is non-singular, so that
$SF(\Delta)^\ast\cong\oplus_{i=1}^u\Z\si_i$, then condition (X) says
that at least one of the (non-negative) integers
$l_z(\si_{i_1}),\dots,l_z(\si_{i_j})$ must be zero whenever
$\si_{i_1},\dots,\si_{i_j}$ do not lie in a single cone.

Second, we may interpret topologically the integer
$\sum_z l_z(\tau)$  as the class $f^\ast [\tau]\in H^2 S^2$. Since
the image of the inclusion $\iota:\Z^r\to SF(\Delta)$ is zero in
$SF(\Delta)/\Z^r\cong H^2X$,  we have:

\noindent(D)\qquad The vector $D=\sum_zl_z\in \SF(\De)^\ast$
belongs to the kernel of the map
\newline  $\iota^\ast:SF(\Delta)^\ast\to(\Z^r)^\ast$.

\noindent It follows from the identification
$\pi_2X\cong\Ker\iota^\ast$ (when $H_2X$ is torsion free; see \S 2)
that we may regard $D$ as the homotopy class of $f$. We shall write
$\Map_D(S^2,X)$ for this component of the space of continuous maps, and
$\Hol_D(S^2,X)$ for its subset consisting of holomorphic maps.

Condition (X) is a local condition, in the sense that it is purely \ll
label-theoretic\rrr. It depends only on the toric variety $X$.
Condition (D), on the other hand, is a global condition, which
depends on $f$.  We shall show that (X) and (D) are the
{\it only} conditions on the configuration
associated to $f$.

\noindent {\bf Definition:} $\QXD(\C)$ denotes the space of
configurations of
distinct points in $\C$ with labels in $\SF(\De)^\ast$ such that
conditions
(X) and (D) are satisfied (for a fixed $D\in SF(\Delta)^\ast$).

\noindent It should be noted that $Q^X_D(\C)$ is in general a topologically
non-trivial subspace of the contractible space
$Q_D(\C;SF(\Delta)^\ast)$, because condition (X) prevents certain
types of collisions.

\proclaim{Proposition 3.1} Let $X$ be a projective toric variety, such
that $H_2X$ is torsion free.  Then the map
$\al^X:\Hol_D^\ast(S^2,X)\to\QXD(\C)$ is a homeomorphism.
\endproclaim

\demo{Proof}
Let $\th:X\to\CP^N$ be an equivariant embedding with
$\th(\ast)=[1;\dots;1]$.  The fact that $\th$ is equivariant means that
it is
induced by a {\it map of fans} $\phi:\De^X\to\De^{\CP^N}$ (\cite{Od1},
Section 1.5), i.e. a $\Z$-linear homomorphism $\phi:\Z^r\to\Z^N$
whose
$\R$-linear extension carries each cone of $\De^X$ into some cone of
$\De^{\CP^N}$. Here, $\De^X$ is the fan of $X$, and $\De^{\CP^N}$ is the
fan of
$\CP^N$. Let $e_1,\dots,e_r$ and $E_1,\dots,E_N$ be
the standard orthonormal bases of $\R^r$ and $\R^N$, respectively. We
denote the
one-dimensional cones of $\De^X$ by $\si_1,\dots,\si_u$, and those
of
$\De^{\CP^N}$ by $\rho_0,\dots,\rho_N$, as usual.

We have an inclusion map
$\th^\prime:\Hol_D^\ast(S^2,X)\to\Hol_d^\ast(S^2,\CP^N)$,
where $\Hol_d^\ast(S^2,\CP^N)$ denotes the space of holomorphic
maps $f$ of some
degree $d$ (depending on $D$) such that $f(\infty)=[1;\dots;1]$.
We also have a map
$\th^{\prime\prime}:\QXD(\C)\to Q^{\CP^N}_d(\C)$,  induced by
$T^\ast:\SF(\De^X)^\ast\to \SF(\De^{\CP^N})^\ast$, where
$T:\SF(\De^{\CP^N})\to \SF(\De^X)$ is given by composition with
$\phi$. The following diagram is commutative:
$$
\CD
\Hol_D^\ast(S^2,X) @>{\th^\prime}>> \Hol_d^\ast(S^2,\CP^N)\\
@V{\al^X}VV    @V{\al^{\CP^N}}VV\\
Q(\C;SF(\Delta^X)^\ast) @>{\th^{\prime\prime}}>>
Q(\C;SF(\Delta^{\C P^N})^\ast)
\endCD
$$
Note that
the map $\al^{\C P^N}$ gives the well known homeomorphism between
$\Hol^\ast_d(S^2,\C P^N)$ and the space of $(N+1)$-tuples of coprime
monic polynomials of degree $d$.

Next we claim that $\th^{\prime\prime}(Q^X_D(\C))\sub
Q^{\C P^N}_d(\C)$. If $\{(z_\al,l_{z_\al})\}_\al\in Q^X_D(\C)$, then we
must check that the configuration $\{(z_\al,T^\ast(l_{z_\al}))\}_\al$
satisfies the two conditions \ll$(\C P^N)$\rr and \ll (d)\rrr. The
first of these is clear from geometrical considerations. The second
amounts to the condition that $\sum_\al T^\ast(l_{z_\al})$ belongs to
the kernel of the map $SF(\Delta^{\C P^N})^\ast \to (\Z^N)^\ast$. That
this is true follows from the fact that  $\sum_\al l_{z_\al}$ belongs
to the kernel of the map $SF(\Delta^{X})^\ast \to (\Z^r)^\ast$.

Since $\th^\prime$ and $\al^{\CP^N}$ are
injective, it follows from the above diagram that $\al^X$ is injective.
To show that $\al^X$ maps surjectively onto $Q^X_D(\C)$, we must
show that a holomorphic map
$f:S^2\to\CP^N$, which has been constructed from a configuration in
the
image of $\th^{\prime\prime}$, actually factors through the
embedding
$\th:X\to\CP^N$.  To do this, we shall need to describe the maps
$\th^\prime,\th^{\prime\prime}$ more explicitly.
We begin with $\th^\prime$. The embedding $\th:X\to\CP^N$ (and
hence the
map $\th^\prime$) is determined by the restriction
$(\C^\ast)^r\to(\C^\ast)^N$ of $\th$ to the corresponding tori
(\cite{Od1}, Theorem 1.13). This is given by
$(z_1,\dots,z_r)\mapsto(z^{m_1},\dots,z^{m_N})$, where
$\phi(x)=\sum_{i=1}^N\lan m_i,x\ran E_i$, and where $z^{m_i}$ means
$z_1^{(m_i)_1}\dots z_r^{(m_i)_r}$. Thus, $X$ may be
described explicitly as the closure in $\CP^N$ of the set of elements
of the form $[1;z^{m_1};\dots;z^{m_N}]$, $z\in(\C^\ast)^r$.
Now we turn to $\th^{\prime\prime}$.  Let $\{(z_\al,l_{z_\al})\}_\al$
be an element of $Q^X_D(\C)$; its image under $\th^{\prime\prime}$ is
$\{(z_\al,T^\ast(l_{z_\al}))\}_\al$. By the remarks above, this
configuration lies in $Q^{\C P^N}_d(\C)$, so it corresponds to
an $(N+1)$-tuple $(p_0,\dots,p_N)$ of coprime monic polynomials of
degree $d$. The exponent of $z-z_\al$ in $p_i(z)$ is
$T^\ast(l_{z_\al})(\rho_i)= l_{z_\al}(T(\rho_i))$. To find the explicit
form of $p_i$, we have to compute $T(\rho_i)$. By definition we have
$T(\rho_i)=\sum_{j=1}^u\lan \rho_i,\phi(v_j)\ran \sigma_j$. Now,
one may verify by direct calculation that
$\lan \rho_i-\rho_0,\sum_{k=1}^N x_kE_k\ran=x_i$,  so we obtain
$$
T(\rho_i)-T(\rho_0)=\sum_{j=1}^u
\lan \rho_i-\rho_0,\phi(v_j)\ran\sigma_j =\sum_{j=1}^u\lan
m_i,v_j\ran\sigma_j.
$$
Hence the exponent of $z-z_\al$ in $p_i(z)p_0(z)^{-1}$ is
$l_{z_\al}(\sum_{j=1}^u\lan m_i,v_j\ran\sigma_j)=
\sum_{k=1}^r(m_i)_k a^\al_k$, where
$a^\al_k=l_{z_\al}(\sum_{j=1}^u(v_j)_k\sigma_j)$. Observe that
$\sum_{j=1}^u(v_j)_k\sigma_j$ belongs to $SF(\Delta)$, since its
value on $v_i$ is just $(v_i)_k$. Hence $a^\al_k$ is an integer, and we
have
%sum changed to Sigma
$$
p_i(z)p_0(z)^{-1}=
\prod_\al(z-z_\al)^{\Sigma_{k=1}^r(m_i)_k a^\al_k } =
q_1(z)^{(m_i)_1}\dots q_r(z)^{(m_i)_r}=q^{m_i},
$$
where $q_k(z)$ denotes the rational function
$\prod_\al(z-z_\al)^{a^\al_k}$. This completes our explicit
determination of the map $\th^{\prime\prime}$.
It follows immediately from this and the earlier description of
$\th^{\prime}$ that the map represented by $(p_0,\dots,p_N)$ factors
through $X$. Hence $\al^X$ maps surjectively onto $Q^X_D(\C)$, as
required.
We have now shown that $\al^X$ is bijective. It is a homeomorphism
because
it is a restriction of $\al^{\CP^N}$, which is a homeomorphism.
\qed\enddemo

\noindent{\it Example 3.2: Complex projective space $\CP^n$.}
With the notation of Example 2.1, we have
$SF(\De)\cong\oplus_{i=0}^n\Z\rho_i\cong SF(\De)^\ast$.
The map $\iota^\ast$ is given by
$$
\iota^\ast:\sum_{i=0}^n x_i\rho_i \longmapsto
(x_1-x_0,\dots,x_n-x_0).
$$
For $D=\sum_{i=0}^n d\rho_i\in\Ker\iota^\ast$, $Q^X_D(\C)$
consists of all configurations such that the labels
$l_z=\sum_{i=0}^n x_i\rho_i$ satisfy the conditions

\noindent (X) \qquad $x_0,\dots,x_n\ge 0$ and
$x_0\dots x_n= 0$

\noindent (D)\qquad $\sum_z x_0=\dots=\sum_z x_n\ (=d)$.

\noindent(Explanation:
A map $f\in\Hol_d^\ast(S^2,\CP^n)$ may be
identified explicitly with an $(n+1)$-tuple $(p_0,\dots,p_n)$ of monic
polynomials of degree $d$ with no common factor. The divisor
$f^{-1}(P_i)$ is given by the roots of $p_i$. Thus, the
labelled configuration associated to $f$ is the set of distinct roots
$z$ of $p_0\dots p_n$, where the label
$l_z=\sum_{i=0}^nx_i\rho_i$ of $z$ indicates that $z$ is a root
of $p_i$ of multiplicity $x_i$.)

\noindent{\it Example 3.3: The Hirzebruch surface $\Sigma_k$.}
{}From Example 2.2,
$SF(\De)\cong\oplus_{i=1}^4\Z\si_i\cong SF(\De)^\ast$, and
the map $\iota^\ast$ is given by
$$
\iota^\ast:\sum_{i=1}^4 x_i\si_i \longmapsto
(x_1-x_3,x_2+kx_3-x_4).
$$
Conditions (X) and (D) are:

\noindent(X)\qquad $x_1,\dots,x_4\ge0$ and
$x_1x_3=0, \ \ x_2x_4=0$

\noindent(D)\qquad $\sum_zx_1=\sum_zx_3,\ \
\sum_zx_2+k\sum_zx_3=\sum_zx_4$.

\noindent(Explanation: From the embedding
$\Sigma_k\sub\C P^1\times\C P^2$ of Example 2.2, we see that a
map $f\in\Hol^\ast(S^2,\Sigma_k)$ may be identified explicitly
with a $5$-tuple of monic polynomials
$((p_4,p_2p_3^k,p_2p_1^k),(p_1,p_3))$,  such that $p_1,p_3$ are
coprime and $p_2,p_4$ are coprime. The divisor
$f^{-1}(X_i)$ is given by the roots of $p_i$. Thus, the
labelled configuration associated to $f$ is the set of distinct roots
$z$ of $p_1p_2p_3 p_4$, where the label
$l_z=\sum_{i=1}^4x_i\si_i$ of $z$ indicates that $z$ is a root of
$p_i$ of multiplicity $x_i$.)

\noindent{\it Example 3.4: The \ll quadric cone\rr $z_2^2=z_1z_3$ in
$\CP^3$.} (The space of rational curves on this variety was considered
in detail in \cite{Gu1}.)
Consider the fan in $\R^2$ given by the three
two-dimensional cones (and all their faces) spanned by the vectors
$v_1=(1,0),v_2=(-1,2),v_3=(0,-1)$.
It can be shown that this fan arises from the quadric cone
$X$ in $\CP^3$ which is defined by the equation $z_2^2=z_1z_3$.
Indeed, this is an example of a weighted projective space (see
Example 2.3), namely $P(1,1,2)$. The torus $(\C^\ast)^2$ acts on $X$
by $(u,v)\cdot[z_0;z_1;z_2;z_3]=[z_0;uvz_1;uz_2;uv^{-1}z_3]$.
We have
$$\gather
\SF(\De)\cong\{h_1\si_1+h_2\si_2+h_3\si_3\ \vert\
h_1,h_2,h_3\in\Z, h_1+h_2\in 2\Z\}\\
\SF(\De)^\ast\cong\{x_1\si_1+x_2\si_2+x_3\si_3\ \vert\
x_1,x_2\in\tfrac12\Z, x_3,x_1+x_2\in\Z\},
\endgather
$$
where we have identified $SF(\Delta)^\ast$ in an obvious way with a
lattice in $\R\sigma_1\oplus\R\sigma_2\oplus\R\sigma_3$.
The map $\iota^\ast$ is given by
$$
\iota^\ast:x_1\si_1+x_2\si_2+x_3\si_3 \longmapsto
(x_1-x_2,2x_2-x_3).
$$
Thus, for $D=d\si_1+d\si_2+2d\si_3$, we see that $Q^X_D(\C)$
consists of all configurations such that the labels
$l_z=x_1\si_1+x_2\si_2+x_3\si_3$ satisfy the conditions

\noindent (X)\qquad $x_1,x_2,x_3\ge0$ and
$x_1x_2x_3=0$

\noindent (D)\qquad $2\sum_z x_1=2\sum_z x_2 = \sum_z x_3\
(=2d)$.

\noindent{\it Example 3.5: The weighted projective space $P(1,2,3)$
(see Example 2.3).} This is a del Pezzo surface.
A suitable fan is the one generated by the vectors
$v_1=(1,0),v_2=(0,1),v_3=(-2,-3)$.  It can be shown (see Appendix 2) that
$$
\SF(\De)\cong\{h_1\si_1+h_2\si_2+h_3\si_3\ \vert\
h_1,h_2,h_3\in\Z, h_2+h_3\in 2\Z,h_1+2h_3,2h_1+h_3\in3\Z\}
$$
$$\align
\SF(\De)^\ast\cong\{x_1\si_1+x_2\si_2+x_3\si_3\ &\vert\
x_1\in\tfrac13\Z,x_2\in\tfrac12\Z,x_3\in\tfrac16\Z,\\
&\ \ \ \ \ \ \ \ x_1+4x_3,2x_1+2x_3,x_2+3x_3,x_1+x_2+x_3\in\Z\},
\endalign
$$
and that the map $\iota^\ast$ is given by
$$
\iota^\ast:x_1\si_1+x_2\si_2+x_3\si_3 \longmapsto
(x_1-2x_3,x_2-3x_3).
$$
For $D=2d\si_1+3d\si_2+d\si_3$, we see that $Q^X_D(\C)$ consists
of all configurations such that the labels
$l_z=x_1\si_1+x_2\si_2+x_3\si_3$ satisfy the conditions

\noindent (X)\qquad $x_1,x_2,x_3\ge0$ and
$x_1x_2x_3=0$

\noindent (D)\qquad $\sum_z x_1=2d$, $\sum_z x_2 = 3d$, $\sum_z
x_3=d$.

\noindent{\it Example 3.6: The tetrahedral complex (see Example
2.5).}
We shall use the following notation:
$$\gather
v_{12}=(1,1,-1)\quad v_{13}=(1,-1,1)\quad v_{23}=(-1,1,1)\quad
v_{123}=(1,1,1)
\\
v_{12}^\prime=(-1,-1,1)\quad v_{13}^\prime=(-1,1,-1)\quad
v_{23}^\prime=(1,-1,-1)\quad v_{123}^\prime=(-1,-1,-1)
\endgather
$$
and we shall write $\si_\ast=\R_{\ge0}v_\ast$,
$\si_\ast^\prime=\R_{\ge0}v_\ast^\prime$, where $\ast$ ranges over
the subscripts $12,13,23,123$. One obtains the identification
$$
SF(\Delta)\cong\{\sum h_\ast\si_\ast+\sum h_\ast^\prime
\si_\ast^\prime\ \vert\  h_\ast,h_\ast^\prime\in\Z,\ (H)\},
$$
where $(H)$ denotes the system of equations
$$\gather
h_{12}+h_{13}=h_{123}+h^\prime_{23},\quad
h_{12}^\prime+h_{13}^\prime=h_{123}^\prime+h_{23}\\
h_{12}+h_{23}=h_{123}+h^\prime_{13},\quad
h_{12}^\prime+h_{23}^\prime=h_{123}^\prime+h_{13}\\
h_{13}+h_{23}=h_{123}+h^\prime_{12},\quad
h_{13}^\prime+h_{23}^\prime=h_{123}^\prime+h_{12}.
\endgather
$$
To identify the dual group $SF(\Delta)^\ast$, let us denote by $V$ the
subspace of $W=(\oplus\R\si_\ast)\oplus(\oplus\R\si_\ast^\prime)$
defined by the equations $(H)$. Then $SF(\Delta)^\ast$ may be
identified with a lattice in the vector space $W^\ast/V^\circ$ (where
$V^\circ=\{f\in W^\ast\ \vert\ f(V)=0\}$), namely the lattice of
functionals which take integer values on $SF(\Delta)$.  Since
$V^\circ$ is generated by the six elements
$\si_{12}+\si_{13}-\si_{123}-\si_{23}^\prime$ etc.,  one may
represent any element of $W^\ast/V^\circ$ by a unique element
$x=\sum x_\ast \si_\ast\in W$. This defines an element of
$SF(\Delta)^\ast$ if and only if it takes integer values on
$SF(\Delta)$. It is easy to check that this is so if and only if
$x_\ast\in\Z$ for all $\ast$. Thus we arrive at the identification
$$
SF(\Delta)^\ast\cong\bigoplus\Z\si_\ast.
$$
Conditions (X) and (D) are as follows:

\noindent (X) \qquad The integers
$x_{12},x_{13},x_{23},x_{12}+x_{13}+x_{123},x_{12}+x_{23}+x_{123},
x_{13}+x_{23}+x_{123}$ are non-negative, but not simultaneously
positive

\noindent (D)\qquad $\sum_z x_{12}= \sum_z x_{13}=\sum_z x_{23}
=-\sum_z x_{123}=d$ where
$D=d\si_{12}+d\si_{13}+d\si_{23}+d\si_{123}$.

\subheading{\S 4  The theorem for non-singular projective toric
varieties}

If $X$ is a non-singular toric variety, then we have
$SF(\Delta)^\ast\cong\oplus_{i=1}^u\Z\si_i$, and so the
homotopy class of a map $S^2\to X$ is given by
an element $D=\sum_{i=1}^ud_i\si_i$ of the kernel of the
homomorphism $\iota^\ast:\oplus_{i=1}^u\Z\si_i\to\Z^r$,
$\si_i\mapsto v_i$, where each $d_i$ is a non-negative integer.
(Recall from \S2 that, in the non-singular case,  $\pi_2X$ may be
identified with $\Ker\iota^\ast$.)  From the identification
$\Hol_D^\ast(S^2,X)\cong Q_D^X(\C)$ of Proposition 3.1,   we have the
following consequences:

\noindent(i) $\Hol_D^\ast(S^2,X)$ is connected

\noindent(ii) the fundamental group of $\Hol_D^\ast(S^2,X)$ is free
abelian of finite rank.

\noindent The first of these follows from the fact that the space
$Q_D^X(\C)$ may be obtained from the affine space of $u$-tuples of
monic polynomials of degrees $d_1,\dots,d_u$, by removing a finite
collection of complex hypersurfaces. The second is proved in the
Appendix of \cite{GKY}.

\proclaim{Theorem 4.1} Let $X$ be a non-singular projective toric
variety. Then the inclusion
$$
\Hol_D^\ast(S^2,X)\to\Map_D^\ast(S^2,X)
$$
is a homotopy equivalence up to dimension $d$, where
$d=\text{min}\{d_1,\dots,d_u\}$ (i.e. this map induces isomorphisms
on homotopy groups in dimensions less than $d$, and an epimorphism
in dimension $d$).\endproclaim

There are two approaches to proving this theorem, depending on how
one views the space $Q^X(\C)$. The first way, which was the original
motivation for this paper, is to view $Q^X(\C)$ as a subspace of
$Q^{\CP^N}(\C)$ via the inclusion map $\th^{\prime\prime}$ (see the
proof of Proposition 3.1). It is the subspace consisting of
configurations of points whose labels satisfy the additional condition
that they belong to the image of the linear transformation $T^\ast$.
This condition arises because of the fact (mentioned earlier) that an
embedding of $X$ in $\CP^N$ may be chosen such that the equations of
$X$ are all of the form \ll monomial $=$ monomial\rrr. It suggests
that the theorem might be proved by interpreting Segal's proof for
$\CP^N$ in terms of labelled configurations, and then imposing the
additional linear condition on the labels. While this is essentially
valid, there are some technical difficulties, and  it is more
convenient (and perhaps more natural) to work with the space
$Q^X(\C)$ directly, without using a particular projective embedding of
$X$.  This is what we shall do.

Let $Q_{d_1,\dots,d_u}^X(U)$ be the space of configurations of
distinct
points $z$ in $U (\sub\C)$ with labels
$l_z\in \SF(\De)^\ast$
which satisfy condition (X), with $\sum_{z}l_z(\si_i)=d_i$ for all
$i$ (where $d_i\ge 0$ for all $i$).  Let
$V$ be an open subset of $\C$ with $U\subset V$.  One may define a
natural
inclusion
$j:Q_{d_1,\dots,d_i,\dots,d_u}^X(U)\to
Q_{d_1,\dots,d_i+1,\dots,d_u}^X(V)$ by adjoining to each
configuration a
fixed point in $V- U$ with the label $\si_i$.

\proclaim{Proposition 4.2}
The inclusion $j:Q_{d_1,\dots,d_i,\dots,d_u}^X(U)\to
Q_{d_1,\dots,d_i+1,\dots,d_u}^X(V)$ is a homotopy equivalence up to
dimension $d_i$.\endproclaim

\demo{Proof} This result is proved in \cite{GKY}. The idea of the proof
is to reduce it to the well known fact that the inclusion $Sp^d(U)\to
Sp^{d+1}(V)$ of symmetric products is a homotopy equivalence up to
dimension $d$.
\qed\enddemo

The space $\QXD(\C)$ of the previous section is of this form.  Let
$D=(d_1,\dots,d_u)$ and $D^\prime=(d^\prime_1,\dots,d^\prime_u)$ be
multi-degrees with $d_i\le d^\prime_i$ for all $i$ (we write $D\le
D^\prime$).  By adjoining a fixed labelled configuration in $V- U$
we obtain an inclusion $j:\QXD(U)\to Q^X_{D^\prime}(V)$.  Evidently
we have:

\proclaim{Corollary 4.3}
The inclusion $j:\QXD(U)\to Q^X_{D^\prime}(V)$ is a homotopy
equivalence up
to dimension $d=\text{min}\{d_1,\dots,d_u\}$.
\qed\endproclaim

We shall introduce a stabilized space using the idea of \cite{Se}.
Let  $\xi=\{(z_i,l_i)\ \vert\ i=1,2,\dots\}$ be a sequence of points of
$Q^X(\C)$ with $z_i\to\infty$.

\noindent{\bf Definition:} $\hat Q_{k_1,\dots,k_u}^{\xi}(\C)$ is the
set of sequences $\{(w_i,m_i)\ \vert\ i=1,2,\dots\}$ of points of
$Q^X(\C)$, which agree with $\xi$ except possibly for a finite number
of terms,  such that $\sum_i(l_i-m_i)_j=k_j$ for $j=1,\dots,u$.  We
shall write $\hat Q_{0}^{\xi}(\C)$ for
$\hat Q_{0,\dots,0}^{\xi}(\C)$.

\noindent Let $D_1<D_2<D_3<\dots$ be a sequence of multi-degrees.
We may
choose open discs $U_1\subset U_2\subset U_3\subset\dots$ in $\C$
and
labelled configurations in each $U_i- U_{i-1}$ so as to obtain a
sequence of inclusions
$$
Q^X_{D_1}(U_1)\longrightarrow Q^X_{D_2}(U_2)\longrightarrow
Q^X_{D_3}(U_3)\longrightarrow\dots
$$
The choice of labelled configurations defines a sequence $\xi$ such
that $\hat Q_0^{\xi}$ is $\bigcup_{i\ge 1}Q^X_{D_i}(U_i)$.
(Since $Q^X_{D_i}(U_i)$ is
homeomorphic to $Q^X_{D_i}(\C)$ and hence to
$\Hol_{D_i}^\ast(S^2,X)$, this
construction may be taken as the definition of the limit
\ll$\lim_{D\to\infty}\Hol_D^\ast(S^2,X)$\rr of \S 1.)

Let $A$ be any subspace of $S^2=\C\cup\infty$, and let $B$ be a
closed
subspace of $A$.  Let $Q^X(A,B)$ denote $Q^X(A)/\sim{}$, where
$Q^X(A)$ is
the space of configurations of distinct points in $A$ with labels in
$\SF(\De)^\ast$ satisfying condition (X), and where two labelled
configurations are defined to be equivalent if they agree on the
complement of $B$ in $A$.
(This is a connected space, if $A$ is connected.) If $\Sigma$ is a
labelled configuration, then
$\Sigma\cap B_z$ defines an element of $Q^X(B_z,\partial B_z)$,
where
$B_z$ is the closed unit disc with centre $z$.  We may identify
$Q^X(B_z,\partial B_z)$ canonically with $Q^X(S^2,\infty)$ and so we
obtain a map
$$
\C\times\QXD(\C)\longrightarrow (Q^X(B_z,\partial
B_z)\longrightarrow)\, Q^X(S^2,\infty),
\quad (z,\Sigma)\longmapsto \Sigma\cap B_z.
$$
  This extends to a continuous map $S^2\times\QXD(\C)\to
Q^X(S^2,\infty)$ with
$(\infty,\Sigma)\mapsto\es$.  The adjoint map
$$
S_D:\QXD(\C)\longrightarrow \Omega^2Q^X(S^2,\infty)
$$
will be called the scanning map.  This is a generalization of a
construction introduced in \cite{Se} for the case $X=\CP^n$. As in
\cite{Se} we
obtain a stabilized map
$$
\hat S:\hat Q_0^{\xi}(\C)\longrightarrow \Omega^2_0Q^X(S^2,\infty),
$$
where $\Omega^2_0$ denotes the component of $\Omega^2$ which
contains the constant maps.

\proclaim{Proposition 4.4}
$\hat S$ is a homotopy equivalence.
\endproclaim

\demo{Proof} This is entirely analogous to the proof in \S 3 of
\cite{Se} for the case $\C P^n$. Another treatment of the same
argument was given in \cite{Gu1} (in the case of the quadric cone) and
in \cite{Gu2} (in the case of $\C P^n$).
\qed
\enddemo

Next we shall examine the relation between $S_D$ and the inclusion
map
$I_D:\Hol_D^\ast(S^2,X)\to\Map_D^\ast(S^2,X)$.

\proclaim{Proposition 4.5} The maps $S_D,I_D$ may be identified
with each other, up to  homotopy.
\endproclaim

\demo{Proof} Let $F^X=\{f:U\to X\ \vert\ f=g\vert_U,\
g\in\Hol(S^2,X)\}$,
where $U$ is the open unit disc in $\C$. The evaluation map
$e:F^X\to X$, $f\mapsto f(0)$ is a homotopy equivalence. Let $\tilde
F^X=\{f\in F^X\ \vert\ f(U)\cap (\TC\cdot\ast)\ne\es\}$. Then $\tilde
F^X$ is obtained from $F^X$ by removing those maps with image in the
complement of $\TC\cdot\ast$ (a subspace of $F^X$ of infinite
codimension), and so the evaluation map $e:\tilde F^X\to X$ is also a
homotopy equivalence.

The action of $\TC$ on $\tilde F^X$ is (by construction) free, in
contrast to
the action of $\TC$ on $X$. (Thus, $\tilde F^X/\TC$ is the
homotopy quotient $X//\TC$.) Let $p:\tilde F^X\to\tilde F^X/\TC$ be
the  natural map.
There is a map $u:\tilde F^X/\TC\to Q^X(\bar U,\partial\bar U)$,
defined by sending the labelled configuration $\Sigma$ which
represents an element $[f]\in\tilde F^X/\TC$ to the labelled
configuration $\Sigma\cap\bar U$.  It may be
shown by an elementary argument as in \cite{Se}, Proposition 4.8 that
$u$ is a homotopy equivalence.

The discussion so far may be summarized in the following diagram:
$$
\CD
\tilde F^X   @>e>>  X\\
@VpVV    @.\\
\tilde F^X/\TC  @>u>>  Q^X(\bar U,\partial\bar U).
\endCD
$$

Consider next the diagram below
$$
\CD
\QXD(\C)  @>{s_D}>>   \Map(\C,\tilde F^X)  @>>>  \Map(\C,X)\\
@V=VV    @VVV   @.\\
\QXD(\C)   @>{\bar s_D}>>  \Map(\C,\tilde F^X/\TC)   @>>>
\Map(\C,Q^X(\bar U,\partial\bar U))
\endCD
$$
where $s_D(f)(z)$ is the map $w\mapsto f(w+z)$ in $\tilde F^X$, for
$f\in\Hol_D^\ast(S^2,X)\cong\QXD(\C)$. (The right hand side of the
diagram
is induced from the previous diagram, and $\bar s_D$ is the map
induced by
$s_D$.)

Observe that $\Map(\C,\ \ \ )$ can be replaced by $\Map^\ast(S^2,\ \ \
)$, i.e.
all the relevant maps extend from $\C$ to $\C\cup\infty=S^2$ (as
based
maps).  Thus we obtain the following commutative diagram, where the
suffix $D$ denotes the appropriate component:
$$
\CD
\QXD(\C)  @>{s_D}>>   \Map_D^\ast(S^2,\tilde F^X)  @>{\Omega^2e}>>
\Map_D^\ast(S^2,X)\\
@V=VV    @V{\Omega^2p}VV   @.\\
\QXD(\C)   @>{\bar s_D}>>  \Map_D^\ast(S^2,\tilde F^X/\TC)
@>{\Omega^2u}>>
\Map_D^\ast(S^2,Q^X(\bar U,\partial\bar U)).
\endCD
$$
The top row is the inclusion map $I_D$, and the bottom row is the
scanning
map $S_D$. The proof of the proposition is completed by noting that
$\Omega^2e,\Omega^2u$ are homotopy equivalences (because $e,u$
are), and
that $\Omega^2p$ is a homotopy equivalence because $p$ is a
fibration with
fibre $(\CS)^r$.
\qed\enddemo

Corollary 4.3, Proposition 4.4 and Proposition 4.5 constitute a proof of
Theorem 4.1.

\subheading{\S 5 Sketch of the theorem in the general case}

In this section we shall sketch how Theorem 4.1 may be
extended to arbitrary compact toric varieties. As it seems
hard to give a single general statement, we shall first obtain a
result under certain special assumptions, and then explain very
briefly how to proceed when the assumptions are not
satisfied. There are two independent parts to the result.
First, we must show:

\noindent (I)\quad The inclusion
$\Hol_D^\ast(S^2,X)\to\Map_D^\ast(S^2,X)$ induces a homotopy
equivalence in the limit $D\to\infty$.

\noindent Then we must find an integer $n(D)$ such that:

\noindent (II)\quad The inclusion
$\Hol_D^\ast(S^2,X)\to\Hol_{D^\prime}^\ast(S^2,X)$
is a homotopy equivalence up to dimension $n(D)$.

\noindent Part (I) can be carried out by making
technical modifications  to the argument of \S 4, as we shall show in
Theorem 5.1. On the other hand, part (II) needs
a new idea, which we give in Theorem 5.2.

To carry out (I), we need to define an appropriate stabilization
procedure. Let $SF(\Delta)^\ast_{\ge0}$  denote the non-negative
elements of $SF(\Delta)^\ast$, i.e. those which take non-negative
values on positive divisors.
Let $D_1, D_2, D_3,\dots$ be a sequence in $SF(\Delta)^\ast_{\ge0}$,
such that each $D_i-D_{i-1}=l_i$ is a valid label; we write
$D_1<D_2<D_3\dots\ $.  As in the non-singular case, we may choose
open discs $U_1\subset U_2\subset U_3\subset\dots$ in $\C$ and a
sequence $\xi=\{(z_i,l_i)\ \vert\ i=1,2,\dots\}$ of labelled points
(with
$z_i\in U_i- U_{i-1}$) so as to obtain a sequence of inclusions
$$
Q^X_{D_1}(U_1)\longrightarrow Q^X_{D_2}(U_2)\longrightarrow
Q^X_{D_3}(U_3)\longrightarrow\ \dots
$$
We obtain a stabilized space
$\hat Q^\xi_0(\C)=\cup_{i\ge0}Q^X_{D_i}(U_i)$ in the usual way. The
inclusion $Q^X_{D_i}(U_i)\to Q^X_{D_{i+1}}(U_{i+1})$ may be regarded
(up to homotopy) as a map $\Hol^\ast_{D_i}(S^2,X)\to
\Hol^\ast_{D_{i+1}}(S^2,X)$, when $D_i,D_{i+1}\in\Ker\iota^\ast$.
We then have:

\proclaim{Theorem 5.1} Let $X$ be a projective toric variety, such
that $H_2X$ is torsion free, and such that the configuration spaces
$Q_D^X(\C)$ are (non-empty and) connected for all $D\in
SF(\Delta)^\ast_{\ge0}$. Then the inclusion
$$
\lim_{D\to\infty}\Hol^\ast_D(S^2,X)\longrightarrow\lim_{D\to\infty}
\Map^\ast_D(S^2,X)\simeq\Map^\ast_0(S^2,X)
$$
is a homotopy equivalence.
\endproclaim

\demo{Proof} This is similar to the proof of Propositions 4.4 and 4.5,
so we shall just point out the new features. For the proof of
Proposition 4.4 one needs to know

\noindent (a) $Q_D^X(\C)$ is (non-empty and) connected for all $D\in
SF(\Delta)^\ast_{\ge0}$,

\noindent (b) the components of $\Map^\ast(S^2,Q^X(S^2,\infty))$ are
indexed by $SF(\Delta)^\ast$, i.e. $\pi_2 Q^X(S^2,\infty)\cong
SF(\Delta)^\ast$, and

\noindent (c) $\pi_1 Q_D^X(\C)$ is abelian for $D$ sufficiently large

\noindent (see the proof of Lemma 3.4 of \cite{Se}).  These three
points are easy to establish when $X$ is smooth, but are not
immediately obvious when $X$ is singular. We have included part (a)
as a hypothesis; in Appendix 1 to this paper we give a method to
determine when this hypothesis is satisfied. Part (b) follows from
the existence of an isomorphism $SF(\Delta)\to H^2 (X//T^{\C})$, which
is given by assigning to a $T^{\C}$-equivariant line bundle on $X$ the
first Chern class of the associated line bundle on $X//T^{\C}$. The
truth of part (c) is explained in Appendix 1. Finally, Proposition 4.5
is proved in exactly the same way, when $X$ is singular.
\qed\enddemo

Our method for (II) depends on the fact that it is possible to choose a
{\it toric resolution} $\th:\hX\to X$. We shall recall briefly this
procedure (see \cite{Fu}, Section 2.7, and \cite{Od1}, Section 1.5).
First, assume that the fan $\Delta$ of $X$ is simplicial, i.e.  for any
cone
$\sigma=\{\sum_{i=1}^ka_iv_i\ \vert\ a_i\ge 0\}$ in $\Delta$, where
$\sigma\cap\Z^r=\Z v_i$, the vectors $v_1,\dots,v_k\in\Z^r$ are
linearly independent.
Let $\si$ be a maximal cone in $\De$. The multiplicity of $\si$ is
defined to be the index of $\bigoplus_{i=1}^k\Z v_i$ in $\Z^r$. By the
criterion for singularity, $X$ is singular if and only if there is a
maximal cone $\si$ in $\De$ of multiplicity greater than one.  In such
a case, there is some $v=\sum_{i=1}^kc_iv_i\in\si\cap\Z^r$ such that
$0\le c_i<1$ for all $i$. Let $\hDe$ be the fan obtained by
sub-dividing $\De$ in the obvious way, i.e. by replacing $\si$ by the
joins of $\R_{\ge0} v$ with all the faces of $\si$. This is a fan
corresponding to a toric variety $\hX$ which is \ll less singular\rr
than $X$. This process of \ll inserting a ray\rr may be repeated
finitely many times, to obtain a non-singular variety $\hX$ and an
equivariant map $\th:\hX\to X$ which is a resolution of $X$.  Finally, if
the fan $\Delta$ is not simplicial, it is easy to see that $\Delta$ may
be made simplicial by inserting suitable rays.

The configuration space for $X$ is related to the configuration space
for $\hX$ by the (set theoretic) formula
$$
Q^X_D(\C)\quad=\bigcup_{\hD\in\th_\ast^{-1}(D)}
Q^{\hX}_{\hD}(\C),
$$
where $\th_\ast:\pi_2\hX\to\pi_2 X$ is the homomorphism induced
by $\th:\hat X\to X$.
The right hand side of this
formula inherits the topology of the left hand side; in Appendix 2
we shall give a more concrete description of this topology.
The idea of our method is to use the
fact that a result of the required type is known for
the spaces $Q^{\hX}_{\hD}(\C)$ (Proposition 4.2, Corollary 4.3).

We shall begin by considering the
special case where  $\De,\hDe$ are simplicial and
$\hDe$ is obtained from $\De$ by inserting a single vector
$v^\prime=\sum_{i=1}^kc_iv_i$ (with $0\le c_i <k$) into a
$k$-dimensional cone $\si$ spanned by vectors $v_1,\dots,v_k$. Let
$v_{k+1},\dots,v_u$ be the generating vectors of the remaining
one-dimensional cones of $\De$. We have
$\SF(\hDe)^\ast\cong\Z\si^\prime\oplus(\bigoplus_{i=1}^u\Z\si_i)$
and we may identify $\SF(\De)^\ast$ with a subspace of
$\bigoplus_{i=1}^u\R\si_i$. The map $T^\ast:\SF(\hDe)^\ast\to
\SF(\De)^\ast$ is given by
$T^\ast(x\si^\prime + \sum_{i=1}^ux_i\si_i)=
\sum_{i=1}^k(x_i+c_ix)\si_i+\sum_{i=k+1}^ux_i\si_i$.
Let $D=\sum_{i=1}^u e_i\si_i\in\Ker\iota^\ast$. From the form of
$T^\ast$, we see that
$$
\th_\ast^{-1}(D)=\{d\si^\prime+\sum_{i=1}^k(e_i-
c_id)\si_i+\sum_{i=k+1}^u e_i\si_i\in\Ker\hat\iota^\ast\ \vert\
0\le d\le \min\{e_i/c_i\ \vert\ c_i\ne 0\}\ \}.
$$
The possible values of $d$ here are of the form
$d_0,d_0+l,\dots,d_0+bl,\dots,d_0+ml$,
for some non-negative integers $d_0,m$.
We denote by $\hD_b$ the element of $\th_\ast^{-1}(D)$
corresponding to $d=d_0+bl$.

{}From now on we shall denote simply by $Q^X_D$ the spaces
$Q^X_D(\C)$ or $Q^X_D(U)$.
For $D,D^\pr\in\Ker\iota$, we write $D^\pr\ge D$ if
$D^\prime=D + \sum_{i=1}^kac_il\si_i  +
\sum_{i=k+1}^ua_i\si_i \in \Ker\iota^\ast$, where
$a,a_{k+1},\dots,a_u$ are non-negative.
We then have a stabilization map $s:Q^X_D\to Q^X_{D^\prime}$.

\proclaim{Theorem 5.2} Let $X$ be a projective toric variety which
admits a resolution $\hX\to X$ of the above form. Let
$D=\sum_{i=1}^u e_i\si_i$ and let $D^\pr\ge D$.
Then the stabilization map $s:Q^X_D\to Q^X_{D^\prime}$ is
a homotopy equivalence up to dimension
$n(D)=\min\{e,e_{k+1},\dots,e_u\}$, where
$$
e=\max\{\min\{d_0+jl,e_1-c_1(d_0+jl),\dots,e_k-c_k(d_0+jl)\}\
\vert\ j=0,1,\dots,m\}.
$$
Moreover, $\lim_{D\to\infty}n(D)=\infty$.
\endproclaim

\noindent
Later we shall give an example to show that the hypothesis on
the resolution is not a serious restriction (Example 5.4).

\demo{Proof} Let
$Q^X_{D;j}=\bigcup_{b=j}^m Q^{\hX}_{\hD_b}$. The stabilization map
$s:Q^X_D\to Q^X_{D^\prime}$ induces a stabilization map
$s_j:Q^X_{D;j}\to Q^X_{D^\prime;j+a}$ for each $j$.
We claim that

\noindent $(\ast)$\quad The stabilization map $s_j$ is a homology
equivalence up to dimension
$n(D;j)=\min\{d_0+jl,e_{k+1},\dots,e_u\}$, and

\noindent $(\ast\ast)$\quad The inclusion
$Q^X_{D;j}\to Q^X_{D;0}=Q^X_D$ is
a homology equivalence up to dimension
$m(D;j)=\min\{e_1-c_1(d_0+jl),\dots,e_k-c_k(d_0+jl)\}$.

\noindent These statements imply that $s$ is a homology equivalence up
to dimension  $\min\{m(D;j),$
\newline
$m(D^\prime,j+a),n(D;j)\}=\min\{m(D;j),n(D;j)\}$.  By choosing
$j$ so as to maximize this number,  we obtain the stated
value of $n(D)$, and it is easy to verify that
$\lim_{D\to\infty}n(D)=\infty$. To prove the theorem, therefore, we must
prove $(\ast)$ and $(\ast\ast)$, and then show that \ll homology\ll
can be replaced by \ll homotopy\rrr.

Statements $(\ast)$ and $(\ast\ast)$ are analogous to
statements (3) and (4) in the proof of Proposition 3.2
of \cite{Gu1}, so we shall just summarize their proofs.

\noindent{\it Proof of $(\ast)$.}
We have $\th_\ast^{-1}(D^\prime) = \{\hD^\prime_0,
\hD^\prime_1,\dots, \hD^\prime_{m+a}\}$, and (for $b=0,1,\dots,m$)
$\hD_{b+a}^\prime=
\hD^\prime_b+al\si^\prime+\sum_{i=k+1}^ua_i\si_i.$
For each $b=j,\dots,m$, we have a stabilization map $\hat
s_b:Q^{\hX}_{\hD_b}\to Q^{\hX}_{\hD_{b+a}^\prime}$, defined by adding
a point of $U^\prime- U$ with label
$\hD_{b+a}^\prime-\hD_b$. The stabilization map $s_j:Q^X_{D;j}\to
Q^X_{D^\prime;j+a}$ is defined by adding the same point with label
$D^\prime-D=\sum_{b=j}^mT^\ast(\hD_{b+a}^\prime-\hD_b)$. These
are compatible, in the sense that the following diagram is
commutative:
%bigcup changed to cup
$$\CD
Q^X_{D;j}  @>=>>  \bigcup_{b=j}^m Q^{\hX}_{\hD_b}\\
@V{s_j}VV  @VV{\cup_{b=j}^m\hat s_b}V\\
Q^X_{D^\prime;j+a}  @>=>>  \bigcup_{b=j}^m Q^{\hX}_{\hD^\prime_{b+a}}
\endCD
$$
By Corollary 4.3, each map $\hat s_b$ is a homology equivalence up to
dimension $\min\{d_0+bl,e_{k+1},\dots,e_u\}$. By the Mayer-Vietoris
argument used in \cite{GKY}, Theorem 2.5, it follows that
$\bigcup_{b=j}^m\hat s_b$ (and hence $s_j$) is a homology
equivalence up to dimension $\min\{d_0+jl,e_{k+1},\dots,e_u\}$.

\noindent{\it Proof of $(\ast\ast)$.}
It suffices to prove that the inclusion  $Q^X_{D;j}\to
Q^X_{D;j-1}$ is a homology equivalence up to dimension $m(D;j)$. To
prove this, we shall use an identification
$$
Q^X_{D;j-1}\simeq Q^X_{D;j}\cup{}_f\ P^{\hX}_{\hD_{j-1}},
$$
where $P^{\hX}_{\hD_{j-1}}$ denotes the configuration space defined
in exactly the same way as $Q^{\hX}_{\hD_{j-1}}$ except that
condition $(\hX)$ is relaxed:  for any point $z$
with label $l$, the integers $\{l(\si_i)\ \vert\ c_i\ne 0\}$ are
allowed to be simultaneously zero.  We define $P^{\hX;1}_{\hD_{j-1}}$
to be the closed subspace of $P^{\hX}_{\hD_{j-1}}$ consisting of
configurations for which the label of at least one point satisfies the
condition $l(\si_i)\ge1$, for all $i$ such that $c_i\ne0$. (It follows
then that $l(\si_i)\ge lc_i$.) We have $Q^{\hX}_{\hD_j}=
P^{\hX}_{D_{j-1}}-P^{\hX;1}_{D_{j-1}}$, and the attaching map $f$ is
the natural map
$$
f:P^{\hX;1}_{\hD_{j-1}}\longrightarrow Q^X_{D;j}.
$$
Hence, the lemma is equivalent to the assertion that the
inclusion $p_{j-1}:
P^{\hX;1}_{\hD_{j-1}}\longrightarrow
P^{\hX;0}_{\hD_{j-1}}=P^{\hX}_{\hD_{j-1}}$ is a homology equivalence
up to dimension $m(D;j)$.
To prove this assertion, we use the stabilization map
$s:P^{\hX;0}_E\longrightarrow P^{\hX;1}_{E+F}$
where $E\in \SF(\hDe)^\ast$ and $F=\sum_{i=1}^kc_il\si_i$.
Consider the composition
$$\CD
P^{\hX;0}_{\hD_{j-1}-F}  @>s>>
P^{\hX;1}_{\hD_{j-1}}  @>{p_{j-1}}>>
P^{\hX;0}_{\hD_{j-1}}.
\endCD
$$
This is homotopic to a stabilization map of the type of Proposition
4.2, hence is a homology equivalence up to dimension
$\min\{e_1-c_1(d_0+(j-1)l)-c_1l,\dots,${}$
e_k-c_k(d_0+(j-1)l)-c_kl\}$.
{}From this we conclude that the map $p_{j-1}$ induces surjections in
homology up to dimension $m(D;j)$.
Next consider the composition
$$\CD
P^{\hX;1}_{\hD_{j-1}}  @>{p_{j-1}}>>
P^{\hX;0}_{\hD_{j-1}}  @>s>>
P^{\hX;1}_{\hD_{j-1}+F}.
\endCD
$$
This is homotopic to the stabilization map
$P^{\hX;1}_{\hD_{j-1}}\to
P^{\hX;1}_{\hD_{j-1}+F}$. By the method of \cite{Gu1}, Proposition
3.2, it may be deduced from Proposition 4.2 that this is a homology
equivalence up to dimension
$\min\{e_1-c_1(d_0+(j-1)l)-c_1l,\dots,${}$
e_k-c_k(d_0+(j-1)l)-c_kl\}.$
{}From this we conclude that the map $p_{j-1}$ induces injections in
homology up to dimension $m(D;j)-1$.
Thus, $p_{j-1}$ is a homology equivalence up to dimension $m(D;j)$,
as required. This completes the proof of $(\ast\ast)$.

To pass from homology to homotopy, we make use of the fact (see
\cite{HH}) that a map induces isomorphisms of homotopy groups if and
only if it induces (a) isomorphisms of homology groups with arbitrary
local coefficients and (b) an isomorphism of fundamental groups. The
stabilization map $Q^X_D\to Q^X_{D^\prime}$ satisfies (a), because
the above argument  for homology with integer coefficients extends
word for word to the case of arbitrary local coefficients: the basic
ingredients were the Mayer-Vietoris exact sequence and the exact
sequence of a pair, together with Proposition 4.2. For (b), we combine
the homology statement with the fact that $\pi_1 Q^X_D$ is abelian
(see Appendix 1).
\qed\enddemo

\noindent{\it Example 5.3: The \ll quadric cone\rr $z_2^2=z_1z_3$ in
$\CP^3$ (see Example 3.4).} Here, $\pi_2 X\cong\Z$.
The variety has one singular point, $[1;0;0;0]$, which
corresponds to the cone $\si$ spanned by $v_1,v_2$. The multiplicity
of $\si$ is $2$. A toric resolution may be obtained by \ll inserting\rr
the vector $v^\pr=\frac12v_1+\frac12v_2=(0,1)$. The corresponding
variety $\hX$ is the Hirzebruch surface $\Sigma_2$ (cf. Example 2.2).

Let $D=g\si_1+g\si_2+2g\si_3\in\Ker\iota^\ast$,
with $g\in\frac12\Z$. For simplicity,
let us assume that $g\in\Z$ (the case where $g-\frac12\in\Z$ is
similar). Then
$\th_\ast^{-1}(D)=\{\hD_0,\hD_1,\dots,\hD_g\}$, where
$\hD_b = (g-b)\si_1 +
2b\si^\pr + (g-b)\si_2 + 2g\si_3$. Here we have $c_1=c_2=\tfrac12$,
$l=2$, $d_0=0$, $m=g$. We have $D^\prime = D + a\si_1 + a\si_2 +
2a\si_3$. By Theorem 5.2, the stabilization map
$Q^X_D\to Q^X_{D^\prime}$ is a homotopy equivalence up to dimension
$\min\{e,2g\}$, where $e=\max\{\min\{2j,g-j\}\ \vert\
j=0,1,\dots,g\} = [\frac23g]$. We conclude that $n(D)=[\frac23g]$ in
this case.

For any projective toric variety $X$, the {\it method} of
Theorem 5.2 may
be used after factoring a resolution $\hX\to X$ into maps of the above
form. Rather than attempt to give a general formula for $n(D)$,
however, we shall just illustrate the method in the
following particular but
non-trivial case.

\noindent{\it Example 5.4: The weighted projective space $P(1,2,3)$
(see Examples 2.3,3.5).} Here, $\pi_2 X\cong\Z$.
The cone spanned by $v_2,v_3$ has multiplicity $2$,
and the cone spanned by $v_1,v_3$ has multiplicity $3$.
By inserting $v^\pr=\frac12v_2+\frac12v_3=(-1,-1)$ we can resolve the
singular point represented by the first cone. By inserting
$v^{\pr\pr}=\frac23v_1+\frac13v_3=(0,-1)$, we replace the second cone
by two maximal cones, one of which represents a singular point,
namely that spanned by $v^{\pr\pr},v_3$. It has multiplicity $2\ (<3)$.
We resolve this singularity by inserting
$v^{\pr\pr\pr}=\frac12v^{\pr\pr}+\frac12v_3=(-1,-2)$.

To use the method of Theorem 5.2, we need to factor our resolution
$\hX\to X$ into a sequence of three simple resolutions.
Let us re-number the
vectors defining the fan $\hDe$ of the resolution $\hX$ as follows:
$v_1=(1,0)$, $v_2=(0,1)$, $v_3=(-1,-1)$,
$v_4=(-2,-3)$, $v_5=(-1,-2)$, $v_6=(0,-1)$.
We shall use the sequence of resolutions
$$
\hX=X_{356}\longrightarrow X_{36} \longrightarrow X_3
\longrightarrow X
$$
where $X_3$ denotes the toric variety whose fan is obtained from the
fan of $X$ by inserting $v_3$, and so on.

\noindent {\it First step: $X_{356}\to X_{36}$}

We shall (temporarily) write $\hX=X_{356},X=X_{36}$.
Let $D=\sum_{i\ne5}e_i\si_i\in \SF(\De)^\ast$. Here
$e_4,e_6\in\tfrac12\Z$ and $e_1,e_2,e_3,e_4+e_6\in\Z$. The map
$T^\ast:\SF(\hDe)^\ast\to \SF(\De)^\ast$ is given by
$T^\ast(\sum_{i=1}^6x_i\si_i)=
x_1\si_1+x_2\si_2+x_3\si_3+(x_4+\tfrac12x_5)\si_4+(x_6+\tfrac1
2x_5)\si_6$.
We have $\th_\ast^{-1}(D)=\{\hD_0,\hD_1,\dots,\hD_m\}$, where
$\hD_b=\sum_{i=1}^6 d_i\si_i$, and $e_1=d_1$, $e_2=d_2$,
$e_3=d_3$, $e_4=d_4+\tfrac12d_5$, $e_6=d_6+\tfrac12d_5$. Let us
assume that $e_4,e_6\in\Z$. Then we may write $d_5=2b$, hence
$\hD_b=e_1\si_1+e_2\si_2+e_3\si_3+(e_4-b)\si_4+
2b\si_5+(e_6-b)\si_6$
for $b=0,1,\dots,m$.
We have $c_4=c_6=\tfrac12$, $l=2$, $d_0=0$, $m=\min\{e_4,e_6\}$.
This is very similar to the situation of Example 5.3. By Theorem 5.2,
the stabilization map $Q^X_D\to Q^X_{D^\prime}$ is a homotopy
equivalence up to dimension
$\min\{[\frac23e_4],[\frac23e_6],e_1,e_2,e_3\}$.

\noindent {\it Second step: $X_{36}\to X_{3}$}

We write $\hX=X_{36},X=X_{3}$.
Let $D=\sum_{i=1}^4e_i\si_i\in \SF(\De)^\ast$. Here
$e_1,e_4\in\tfrac13\Z$ and $e_2,e_3,e_1+e_4\in\Z$. The map
$T^\ast:\SF(\hDe)^\ast\to \SF(\De)^\ast$ is given by
$T^\ast(\sum_{i\ne 5}x_i\si_i)=${}
$(x_1+\tfrac23x_6)\si_1+x_2\si_2+x_3\si_3+
(x_4+\tfrac13x_6)\si_4$.
Let us assume that $e_1,e_4\in\Z$. Then we have
$\th_\ast^{-1}(D)=\{\hD_0,\hD_1,\dots,\hD_m\}$, where
$\hD_b=(e_1-2b)\si_1+e_2\si_2+e_3\si_3+(e_4-b)\si_4+
3b\si_6$.
In this situation, $c_1=\frac23$, $c_4=\frac13$, $l=3$, $d_0=0$,
$m=\min\{[\frac12e_1],e_4\}$.
We may now apply (the method of) Theorem 5.2.
A slight strengthening of the result of the first
step is needed here, namely that the \ll individual\rr
stabilization map
$s_i:Q^{\hX}_D\to Q^{\hX}_{D+\si_i}$
is a homotopy equivalence up to dimension $e_i$ for
$i=1,2,3$, and up to dimension $[\frac23e_i]$ for $i=4,6$. (We
omit the proof of this result.) We find that the
stabilization map $Q^X_D\to Q^X_{D^\prime}$ is a homotopy
equivalence up to dimension $\min\{e,e_2,e_3\}$, where
$$
e=\max\{\min\{\tfrac23(3j),e_1-2j,[\tfrac23(e_4-j)]\}\ \vert\
j=0,1,\dots,m\} = \min\{[\tfrac12e_1],[\tfrac12e_4]\}.
$$

\noindent {\it Third step: $X_{3}\to X$}

We write $\hX=X_3$.
Let $D=e_1\si_1+e_2\si_2+e_4\si_4\in \SF(\De)^\ast$.  The
conditions on $e_1,e_2,e_4$ were given earlier. The map
$T^\ast:\SF(\hDe)^\ast\to \SF(\De)^\ast$ is given by
$T^\ast(\sum_{i=1}^4x_i\si_i)=${}
$x_1\si_1+(x_2+\tfrac12x_3)\si_2+(x_4+\tfrac12x_3)\si_4$.
Let us assume that $e_2,e_4\in\Z$. Then
$\th_\ast^{-1}(D)=\{\hD_0,\hD_1,\dots,\hD_m\}$, where
$\hD_b=e_1\si_1+(e_2-b)\si_2+2b\si_3+(e_4-b)\si_4$.
This time we have $c_2=c_4=\frac12$, $l=2$, $d_0=0$,
$m=\min\{e_2,e_4\}$.
By the above method we find that the stabilization map
$Q^X_D\to Q^X_{D^\prime}$ is a homotopy equivalence up to dimension
$\min\{e,[\frac12e_1]\}$, where
$$
e=\max\{\min\{2j,e_2-j,[\tfrac12(e_4-j)]\}\ \vert\ j=0,1,\dots,m\}
= \min\{[\tfrac23e_2],[\tfrac25e_4]\}.
$$
In conclusion, we have shown that for $X=P(1,2,3)$ the inclusion
$\Hol_D^\ast(S^2,X)\to\Map_D^\ast(S^2,X)$ is a homotopy equivalence
up to dimension
$n(D)=\min\{[\tfrac12e_1],[\tfrac23e_2],[\tfrac25e_4]\}$, where
$D=e_1\si_1+e_2\si_2+e_4\si_4$. Thus, if
$D=2d\si_1+3d\si_2+d\si_4$, then $n(D)=[\tfrac25d]$.

Finally, we shall indicate how to proceed in the case of a
(compact) toric variety which is not covered by our methods
up to this point. There are two problems to deal with, namely
(i) the description of the connected components of
$\Hol(S^2,X)$ when $H_2X$ has torsion, and (ii) the extension of
all our previous results to the case of a non-projective
toric variety.

Regarding (i), let us consider a (singular, projective, compact)
toric variety $X$, and let us choose $D\in \Ker\iota^\ast$.
There is a surjection
$\de:H_2X(\cong\pi_2X)\to \Ker\iota^\ast(\cong (H^2X)^\ast)$,
whose kernel is the torsion subgroup of $H_2X$.
We define
$$
\Hol_D^\ast(S^2,X)=\{f\in\Hol^\ast(S^2,X)\ \vert\
\de[f]=D\}
$$
where $[f]\in\pi_2X$ is the homotopy class of $f$. The space
$Q_D^X(\C)$ is defined in the usual way, and we have
$\Hol_D^\ast(S^2,X)\cong Q_D^X(\C)$. The problem now is to
describe the connected components of $Q_D^X(\C)$.
Although this may be done in any particular case
by the method of Appendix 1,
we are unable to give a general statement. Therefore, we
shall content ourselves by giving the following example.

\noindent{\it Example 5.5: The tetrahedral complex (see Examples
2.5, 3.6).} In this case, the group $H_2X$ has been
computed by M. McConnell (in a private
communication) to be $\Z\oplus(\Z/2)^{11}$. For
$D=d\si_{12}+d\si_{13}+d\si_{23}+d\si_{123}\in\Ker\iota^\ast$, the
method described in Appendix 1 shows that $Q^X_D$
has eight components, if $d\ge3$. Four of the components correspond
to holomorphic maps whose image lies in one of four copies of $\C
P^1$, each of which is given by conditions of the form $z_i=z_j=z_k$
with $i\in\{0,1\}$, $j\in\{2,3\}$, $k\in\{4,5\}$.  The other four
components consist of \ll full\rr holomorphic maps.

In the terminology of Appendix 1, the \ll simple\rr labels here are:
$$
\gather
l_{12}=(1,1,0,-1),\quad l_{13}=(1,0,1,-1),\quad
l_{23}=(0,1,1,-1),\quad l_{123}=(1,1,1,-2)\\
l_{12}^\prime=(0,0,1,0),\quad l_{13}^\prime=(0,1,0,0),\quad
l_{23}^\prime=(1,0,0,0),\quad l_{123}^\prime=(0,0,0,1).
\endgather
$$
To define a stabilization procedure, let us consider the sequence of
labels $l_1,l_2,l_3,\dots$, where $l_1,\dots,l_8$ are the above
simple labels and $l_i=l_{i-8}$ for $i>8$. Let $z_1,z_2,z_3,\dots$ be
a sequence of points such that $z_i\in U_i- U_{i-1}$, as above. With
$D_i=l_1+\dots+l_i$, we define $\underline{Q}_{\,D_i}^X(U_i)$ to be
the component of  $Q_{D_i}^X(U_i)$ which contains the configuration
$\{(z_j,l_j)\ \vert\ j=1,\dots,i\}$. Then for $i=8j$, the space
$\underline{Q}_{\,D_i}^X(U_i)$ may be identified up to homotopy with
a distinguished component $\underline{\Hol}_{\,4j}^\ast(S^2,X)$ of
$\Hol_{4j}^\ast(S^2,X)$. Let $\underline{\Map}_{\,4j}^\ast(S^2,X)$
denote the component of $\Map(S^2,X)$ containing the image of
$\underline{\Hol}_{\,4j}^\ast(S^2,X)$.
The method of the proof of Theorem 5.1 then gives a homotopy
equivalence
$\lim_{j\to\infty}\underline{\Hol}_{\,4j}^\ast(S^2,X)
\longrightarrow\lim_{j\to\infty}\underline{\Map}_{\,4j}^\ast(S^2,X)$.
Thus, we obtain a {\it modified} version of our main theorem
in this case.

Regarding (ii), the main question is whether Proposition 3.1
can be proved in the case of a non-projective toric variety.
It suffices to consider the non-singular case,
because singular varieties
may be dealt with by using a resolution.
A proof may be obtained from the following construction of a
toric variety $X$ from its fan $\De$, described in
\cite{Co1},\cite{Au}. Let
$$
Z=\{\sum x_i\si_i\in\bigoplus_{i=1}^u\C\si_i\ \vert\ \prod_{i\in
I_\si^c}x_i=0\ \text{for all}\ \si\in\De\}
$$
where $I_\si=\{i\ \vert\ \si_i\subseteq\si\}$ and $I_{\si}^c$ is the
complement of $I_{\si}$ in $\{1,\dots,u\}$. Let $G$ be the kernel of
the map
$$
(\C^\ast)^u\cong\bigoplus_{i=1}^u\C\si_i/
\bigoplus_{i=1}^u\Z\si_i\longrightarrow\C^r/\Z^r
$$
induced by $\iota\otimes\C$. This is an algebraic group, which acts
naturally on  $\bigoplus_{i=1}^u\C\si_i$.

\proclaim{Theorem 5.6 (\cite{Co1},\cite{Au})} Let $X$ be a simplicial
toric variety. The action of $G$ preserves
$\bigoplus_{i=1}^u\C\si_i-Z$, and the quotient space is isomorphic to
$X$. If $X$ is
non-singular (hence, in particular, simplicial), the action of $G$ is
free.
\qed\endproclaim

\noindent This is analogous to the usual description of $\CP^n$ as the
quotient $(\C^{n+1}-\{0\})/(\C^\ast)^n$, which is a special case.

We may now give an alternative proof of Proposition 3.1,
in the non-singular case.
The idea of the proof is that, just as a (based) holomomorphic map
$S^2\to\CP^n$ may be represented by an
$(n+1)$-tuple of monic polynomials which have no common factor, a
(based) holomorphic map $S^2\to X$ may be represented by a
$u$-tuple of monic polynomials $(p_1,\dots,p_u)$ such that
$(p_1(z),\dots,p_u(z))\notin Z$ for all $z\in\C$. Let $I_z=\{i\ \vert\
p_i(z)=0\}$. Now, $(p_1(z),\dots,p_u(z))\notin Z$ if and only if there
is a cone $\si\in\De$ such that $I_z\cap(I_\si)^c=\emptyset$, i.e.
$I_z\subseteq I_\si$. So the condition on $p_1,\dots,p_u$ is that, if
$p_{i_1},\dots,p_{i_j}$ have a common factor, then
$\si_{i_1},\dots,\si_{i_j}$ belong to a single cone of the fan. This is
precisely condition (X) of \S 3. A complete proof, in a more general
context, has been given recently by Cox in \cite{Co2}.

\subheading{APPENDIX 1: $\pi_0Q_D^X$ and $\pi_1Q_D^X$}

For a toric variety $X$ we have given a correspondence between
(based) holomorphic maps $S^2\to X$ and configurations of points
with labels in $SF(\Delta)^\ast$. These labels satisfy condition (X) of
\S 3. The set of all such labels, being a subset of the monoid
$SF(\Delta)^\ast_{\ge0}$, has the structure of a partial monoid,
which we shall denote by $M_X$. In this appendix we shall indicate
briefly how the algebraic structure of $M_X$ determines
$\pi_i\Hol_D^\ast(S^2,X)$ (i.e. $\pi_i Q_D^X(\C)$) for $i=0,1$.

Let $l_1,\dots,l_p$ be the simple elements of $M_X$ (an element is
said to be simple if it cannot be written in the form $m_1+m_2$,
with $m_1$ and $m_2$ both non-zero). To investigate whether
$Q_D^X$ is connected, we note first that any element of $Q_D^X$ may
be moved continuously to a configuration of points whose labels are
all simple. If that configuration contains certain points
$z_1,\dots,z_q$ with labels $l_{i_1},\dots,l_{i_q}$, such that the
sum $l_{i_1}+\dots+l_{i_q}$ is defined (in $M_X$), then it may be
moved continuously to a configuration in which $z_1,\dots,z_q$ are
replaced by a single point $z$ (with label $l_{i_1}+\dots+l_{i_q}$).
Our aim now is to repeat this reduction process until we arrive at a
canonical configuration; if this is possible, then we will have shown
that $Q_D^X$ is
path-connected.

Let us suppose that there exist {\it linearly independent} labels
$m_1,\dots,m_r\in M_X$, such that the above reduction process
eventually leads to a configuration of the form $\{(z_i,k_im_i)\
\vert\ i=1,\dots,r\}$, where $k_1,\dots,k_r$ are non-negative
integers.  Since $D=\sum k_im_i$ (and $m_1,\dots,m_r$ are linearly
independent), the integers $k_1,\dots,k_r$ are determined by $D$.
Therefore,  we have succeeded in moving to a canonical configuration,
and so $Q_D^X$ is connected.

In the examples occurring in this paper it is straightforward to
determine whether $m_1,\dots,m_r$ exist. For example, in the case
of the quadric cone (Examples 3.4,5.3), the simple labels are
$\si_1,\frac12(\si_1+\si_2),\si_2,\si_3$. If
$D=d_1\si_1+d_2\si_2+d_3\si_3$, a suitable choice of
$m_1,\dots,m_r$ would be $\si_1,\frac12(\si_1+\si_2),\si_3$ (if
$d_1\ge d_2$), or $\frac12(\si_1+\si_2),\si_2,\si_3$ (if $d_1\le
d_2$).

The simple labels are also the main ingredient in the computation of
the fundamental group of $Q^X_D$. By a slight generalization of the
argument used in the Appendix of \cite{GKY}, it follows that the
fundamental group is abelian. Moreover, there is one generator for
each pair of simple labels $l_i,l_j$ such that the sum $l_i+l_j$ is
{\it not} defined (in $M_X$). The order of such a generator is the
least positive integer $n$ such that $n(l_i+l_j)\in M_X$.

\subheading{APPENDIX 2: Representation of holomorphic maps
by polynomials}

In the case of a non-singular toric variety
$X$, Proposition 3.1 (and Theorem 5.6) gives a description
of any $f\in\Hol_{D}^\ast(S^2,X)$
as a sequence $(p_1,\dots,p_{u})$
of monic polynomials, where

\noindent (X)\quad $p_{i_1},\dots,p_{i_j}$ are coprime  if
$X_{i_1}\cap\dots\cap X_{i_j}=\es$, and

\noindent (D)\quad $\deg p_i=d_i$, where $D=
\sum_{i=1}^{u}d_i\si_i$.

\noindent The roots of the polynomial $p_i$ represent the
divisor $f^{-1}(\si_i)$. This description is canonical;
it does not depend on any
embedding in projective space. Such an embedding merely converts the
above polynomial description into a more complicated one,
as we have seen in Example 3.3.

In the case of a singular variety $X$,  the divisors
$\si_1,\dots,\si_u$ are not necessarily Cartier divisors, so we
cannot expect the same procedure to work. Instead, let us {\it choose}
a generating set $\tau_1\dots,\tau_v$ for the positive divisors in
$\SF(\De)$, and then define monic polynomials $q_1,\dots,q_v$ by
taking the roots of $q_i$ to represent the divisor
$f^{-1}(\tau_i)$. This is the same as $\hat f^{-1}(T(\tau_i))$,
where $T:SF(\De)\to SF(\hat\De)$ is the map induced by a toric
resolution $\hat X\to X$, and where $\hat f:S^2\to\hat X$ corresponds
to $f:S^2\to X$.
If $T(\tau_i)=\sum_{j=1}^{\hat u}b_{ij}\hat\si_j$, then
we obtain
$$
(q_1,\dots,q_v)=(p_1^{b_{11}}\dots p_{\hat u}^{b_{1\hat u}},
\ \dots\ ,p_1^{b_{v1}}\dots p_{\hat u}^{b_{v\hat u}})=
(p^{b_1},\dots,p^{b_v}).
$$
The proof of Proposition 3.1 shows that
elements of $\Hol_D^\ast(S^2,X)$ may be
identified with $v$-tuples $(p^{b_1},\dots,p^{b_v})$, where the
polynomials $p_1,\dots,p_{\hat u}$ satisfy

\noindent (X)\quad $p_{i_1},\dots,p_{i_j}$ are coprime  if
$\hX_{i_1}\cap\dots\cap\hX_{i_j}=\es$, and

\noindent (D)\quad $\deg p^{b_i}=D(\tau_i)$, $i=1,\dots,v$.

\noindent This illustrates how $\Hol^\ast_D(S^2,X)$ is topologized as
the union of those $\Hol^\ast_{\hD}(S^2,\hX)$ for which
$\th_\ast(\hD)=D$: each $\Hol^\ast_{\hD}(S^2,\hX)$ has its usual
topology, but a collection of roots of polynomials
$p_{i_1},\dots,p_{i_j}$ may coalesce to give a root of another
polynomial $p_{i_k}$,  where $p_{i_k}$ is associated with a ray
which sub-divides the cone associated to $p_{i_1},\dots,p_{i_j}$.

This kind of polynomial description of elements of $\Hol^\ast_D(S^2,X)$
also arises if we consider a suitable embedding of $X$ in
projective space.

\noindent {\it Example A2.1: The quadric cone (see Examples 3.4,5.3).}
Let us choose the
generators $\si_4,2\si_1,\si_1+\si_3,2\si_3$ of $\SF(\De)$. We
have
$T(\si_4)=\si_4$, $T(2\si_1)=2\si_1+\si_2$,
$T(\si_1+\si_3)=\si_1+\si_2+\si_3$, $T(2\si_3)=2\si_3+\si_2$,
so $f\in\Hol_D^\ast(S^2,X)$ is represented by a $4$-tuple of
polynomials
$$
(q_1,q_2,q_3,q_4)=(p_4,p_1^2p_2,p_1p_2p_3,p_2p_3^2)
$$
where

\noindent (X)\quad $p_1,p_3$ are coprime, $p_2,p_4$ are coprime,
and

\noindent (D)\quad $\deg p_4=\deg p_1^2p_2=\deg p_1p_2p_3=\deg
p_2p_3^2=2g$.

\noindent This is in fact the polynomial representation which arises
from the given embedding in $\CP^3$. In other words, as one readily
verifies,  a $4$-tuple $(q_1,q_2,q_3,q_4)$ of coprime monic
polynomials of degree $2g$ satisfies the equation $q_2^2=q_1q_3$ if
and only if it is of the above form.

\noindent {\it Example A2.2: The weighted projective space $P(1,2,3)$
(see Examples 3.5,5.4).}
Let us choose the  generators
$6\si_4$, $\si_1+4\si_4$, $2\si_1+2\si_4$,
$3\si_1$, $\si_2+3\si_4$, $2\si_2$, $\si_1+\si_2+\si_4$
of $\SF(\De)$. (The reason for this choice will become clear in a
moment.) Applying the map $T$, we find that any
$f\in\Hol_D^\ast(S^2,X)$ is represented by a $7$-tuple of polynomials
$(q_1,q_2,q_3,q_4,q_5,q_6,q_7)$ of the form
$$
(\ p_3^3p_4^6p_5^4p_6^2,
\ p_1p_3^2p_4^4p_5^3p_6^2,
\ p_1^2p_3p_4^2p_5^2p_6^2,
\ p_1^3p_5p_6^2,
\ p_2p_3^2p_4^3p_5^2p_6,
\ p_2^2p_3,
\ p_1p_2p_3p_4p_5p_6\ )
$$
where $p_1,p_2,p_3,p_4,p_5,p_6$ satisfy the conditions

\noindent (X)\quad $p_i,p_j$ are coprime except possibly when
$\vert i-j\vert=1$ or $\{i,j\}=\{1,6\}$, and

\noindent (D)\quad $\deg q_i=6g,\ i=1,\dots,7$.

\noindent In \cite{Ha}, Example 10.27, it is shown that $P(1,2,3)$ may
be embedded in $\CP^6$ via the equations
$z_0z_2=z_1^2,\ z_2z_5=z_6^2,\ z_1z_3=z_2^2,\ z_1z_5=z_4z_6$.
Our polynomial representation was chosen to be compatible
with this embedding.

\enddocument
\end